\documentclass[11pt]{amsart}

\usepackage{ae}
\usepackage[alphabetic]{amsrefs}
\usepackage{mathrsfs}

\newcommand{\HypL}{\textbf{(L)}\ }

\renewcommand{\setminus}{\ \rule[2.5pt]{7pt}{1pt}\ }
\newcommand{\ideal}{\lhd}

\newcommand{\glie}{\mathfrak{g}}

\newcommand{\clie}{\mathfrak{c}}
\newcommand{\sllie}{\mathfrak{sl}}
\newcommand{\gl}{\mathfrak{gl}}
\newcommand{\plie}{\mathfrak{p}}
\newcommand{\ulie}{\mathfrak{u}}

\newcommand{\Lie}{\operatorname{Lie}}
\newcommand{\G}{\mathbf{G}}
\newcommand{\Gm}{{\G_m}}
\newcommand{\Aff}{\mathbf{A}}
\newcommand{\Proj}{\mathbf{P}}
\newcommand{\Ad}{\operatorname{Ad}}
\newcommand{\ad}{\operatorname{ad}}

\newcommand{\A}{\mathcal{A}}

\newcommand{\U}{\mathcal{U}}
\newcommand{\V}{\mathcal{V}}
\newcommand{\W}{\mathcal{W}}

\newcommand{\OO}{\mathscr{O}}
\newcommand{\Orbit}{\mathcal{O}}
\newcommand{\MM}{\mathcal{M}}
\newcommand{\NN}{\mathcal{N}}
\newcommand{\LL}{\mathcal{L}}
\newcommand{\VV}{\mathcal{V}}

\newcommand{\B}{\mathcal{B}}

\newcommand{\Z}{\mathbf{Z}}

\newcommand{\mm}{\mathfrak{m}}

\newcommand{\SL}{\operatorname{SL}}
\newcommand{\GL}{\operatorname{GL}}
\newcommand{\SP}{\operatorname{Sp}}
\newcommand{\SO}{\operatorname{SO}}
\newcommand{\Spec}{\operatorname{Spec}}
\newcommand{\tensor}{\otimes}

\newtheorem{prop}{Proposition}
\newtheorem{lem}[prop]{Lemma}
\newtheorem{theorem}[prop]{Theorem}

\theoremstyle{remark}
\newtheorem{example}[prop]{Example}
\newtheorem{question}[prop]{Question}
\newtheorem{rem}[prop]{Remark}

\begin{document}

\title{On the centralizer of the sum of commuting nilpotent elements}
\author{George McNinch}
\address{Department of Mathematics \\ 
         Tufts University \\ 
         503 Boston Avenue \\ 
         Medford, MA 02155 \\ 
         USA}
\email{george.mcninch@tufts.edu}
\date{April 12, 2005}
\thanks{Research of the author supported in part by the US National
  Science Foundation through DMS-0437482}
\dedicatory{To Eric Friedlander, on his 60th birthday.}


\begin{abstract}
  Let $X$ and $Y$ be commuting nilpotent $K$-endomorphisms of a vector
  space $V$, where $K$ is a field of characteristic $p \ge 0$. If
  $F=K(t)$ is the field of rational functions on the projective line
  $\mathbf{P}^1_{/K}$, consider the $K(t)$-endomorphism $A=X+tY$ of
  $V$. If $p=0$, or if $A^{p-1}=0$, we show here that $X$ and $Y$ are
  tangent to the unipotent radical of the centralizer of $A$ in
  $\GL(V)$.  For all geometric points $(a:b)$ of a suitable open
  subset of $\mathbf{P}^1$, it follows that $X$ and $Y$ are tangent to
  the unipotent radical of the centralizer of $aX+bY$. This answers a
  question of J. Pevtsova.
\end{abstract}

\maketitle

Let $G$ be a connected and reductive algebraic group defined over an
arbitrary field $K$ of characteristic $p \ge 0$.  Write $\glie =
\Lie(G)$, and consider the extension field $F=K(t)$ with $t$
transcendental over $K$. For convenience, we fix an algebraically
closed field $k$ containing both $K$ and $t$.

If $X,Y \in \glie(K)$ are nilpotent and $[X,Y]=0$, then $A=X + tY \in
\glie(F)$ is again nilpotent.  Write $C$ for the centralizer of $A$ in
$G$, and write $R_uC$ for the unipotent radical of $C$. Under
favorable restrictions on the characteristic, the groups $C$ and
$R_uC$ are defined over $K(t)$. In this note, I want to answer -- at
least in part -- a question put to me by Julia Pevtsova at the July
2004 meeting in Snowbird, Utah.  With notation as before, this
question may be stated as follows:

\begin{question}
  When is it true that $X,Y \in \Lie R_uC$?
\end{question}

To begin the investigation, the first section of the paper includes some
elementary results concerning $G$-varieties in case the algebraic
group $G$ acts with a finite number of orbits. For the most part, the
use of these results could be avoided in the present application, but
there is perhaps some interest in recording them.

After these preliminaries, I am mainly going to investigate Question 1
in case the $K$-group is $G=\GL(V)$, where $V$ is a finite dimensional
$k$-vector space defined over $K$; this means there is a given
$K$-subspace $V(K)$ for which the inclusion induces an isomorphism
$V(K) \tensor_K k \simeq V$.

The second section contains well-known material on nilpotent orbits,
mainly for the group $\GL(V)$; this material is used in section three
where we prove our main result -- Theorem \ref{theorem:GL(V)} --
giving a partial answer to Question 1 when $G$ is the group $\GL(V)$.
A final section contains some remarks about more general semisimple
groups.

Let me make a few remarks about possible reasons for interest in the
main result of this paper. Pevtsova's interest concerns finite group
schemes over a field $K$ of characteristic $p>0$; see e.g.  \cite{FP}.
Basic but important examples are the commutative, \'etale, unipotent
group schemes; consider e.g.  a constant finite group scheme $E$
which ``is'' an elementary Abelian $p$-group.  If $(\rho,M)$ is a
$K$-representation of $E$, the matrices $1-\rho(g) = \rho(1-g) \in
\operatorname{End}_K(M)$ are nilpotent for $g \in E$. More generally,
if $x$ is in the augmentation ideal of the group algebra $KE$, then
$\rho(x)$ is nilpotent, and is a linear combination of commuting
nilpotent matrices $\rho(1-g)$ for various elements $1 \ne g$ of $E$.
Pevtsova's question was aimed at understanding properties of the
Jordan block structure of suitably generic such $x$.

In a somewhat different direction, if $G$ is a reductive group over
$K$ and $\NN \subset \glie$ denotes the variety of nilpotent
elements, one is interested in studying the subvariety $\V_2 \subset \NN
\times  \NN$ of commuting pairs:
\begin{equation*}
  \V_2 = \{(X_1,X_2) \in \NN^2 \mid [X_1,X_2]=0\};
\end{equation*}
see e.g. \cite{premet-commute}.  Any $K$-point
\begin{equation*}
  x = (X_1,X_2) \in \V_2(K)
\end{equation*}
determines a nilpotent element $A=X_1 + tX_2 \in \glie(F)$ with
$F=K(t)$ as before.  One might hope to exploit the results of this
paper to study properties of the variety $\V_2$.

\section{Groups acting with finitely many orbits}
In this section, we work ``geometrically'' -- i.e. over the
algebraically closed field $k$. The results recorded here are
elementary and without doubt are well-known; however, I don't know of an
adequate reference.

Let $W$ be an irreducible affine $k$-variety with coordinate algebra
$A = k[W]$. [I will identify $k$ varieties with their $k$-points: $W =
W(k)$.]  For an extension field $k'$ of $k$, 
write $W(k')$ for the $k'$-points of $W$, and write $W_{/k'}$ for the
$k'$-variety obtained by extension of scalars:
\begin{equation*}
  W_{/k'} = W \times_{\Spec(k)} \Spec(k').
\end{equation*}

We will be concerned here with the case where a $k$-group acts on
$W$ with a finite number of orbits.

\subsection{Invariance of the number of orbits}

Begin with the following:
\begin{lem}
  \label{lem:open}
  If $W$ is the union $W = W_1 \cup W_2 \cup \cdots \cup W_n$ of
  locally closed subvarieties $W_j$, then $W_i$ is a non-empty open
  subset of $W$ for some $1 \le i \le n$.
\end{lem}

\begin{proof}
  For $1 \le j \le n$, write $W_j = C_j \cap U_j$ where $C_j \subset
  W$ is closed and $U_j \subset W$ is open. Since $W$ is contained in
  the union of the $C_j$ and irreducible, we find $W \subseteq C_i$ for
  some $i$ and the lemma follows.
\end{proof}

Let $G$ be a connected linear algebraic group over $k$ acting by
$k$-automorphisms on the variety $W$.  Let $x \in W = W(k)$, and let
$\Orbit = G.x$.
Since $\Orbit$ is a $k$-variety, one may speak of its
$k'$-points $\Orbit(k') \subseteq W(k')$. On the other hand, one may
regard $x$ as an element of $W(k')$ and form its $G(k')$-orbit. 
\begin{lem} 
  \label{lem:orbit-lemma}
  Let $x$ be as above, and suppose that the extension field $k'$ of
  $k$ is itself algebraically closed. Then we have
  \label{lem:orbit-claim}
  \begin{equation}
    \label{eq:orbit-claim}
      G(k')x = \Orbit(k').
  \end{equation}
\end{lem}

\begin{proof}
  Since $\Orbit$ is locally closed, we may replace $W$ by the closure
  of $\Orbit$, and so suppose $\Orbit$ to be open in $W$. Since $W
  \setminus \Orbit$ is a union of $G$-orbits each of dimension $<\dim
  \Orbit$, $W(k') \setminus \Orbit(k')$ is $G(k')$-stable, and so
  $\Orbit(k')$ is $G(k')$-stable  
  \footnote{One may avoid arguing the $G(k')$-stability of $\Orbit(k')$ by
  applying \cite{DG}*{II.5.3.2(a)}.}
  .  Since $x \in \Orbit(k')$, the
  containment $\subseteq$ of \eqref{eq:orbit-claim} is immediate.
  
  One finds e.g. in \cite{springer-LAG}*{Proposition 1.9.4 and Theorem
    1.9.5} the elementary proof -- which goes back to Chevalley and
  Weil -- that the image $\phi(X)$ of a dominant morphism of affine
  $k$-varieties $\phi:X \to Y$ contains a non-empty open subset of
  $Y$.  That proof shows more precisely that there is some regular
  function $0 \not = f \in k[Y]$ such that $D(f)(k') \subseteq
  \phi(X(k'))$ for each algebraically closed field $k'$ containing
  $k$; here $D(f)$ is the ``distinguished open'' subset of $Y$
  determined by the non-vanishing of $f$.
  
  Apply this now to the (dominant) orbit map $(g \mapsto gx):G \to W$
  to find $0 \not = f \in A=k[W]$ such that $D(f) \subseteq \Orbit$
  and $D(f)(k') \subseteq G(k')x$. Since $\Orbit=Gx$ is a Noetherian
  space, there are elements $g_1,\dots,g_n \in G = G(k)$ such that
  $Gx$ is the union of the $g_iD(f)$.  Then $\Orbit(k')$ is the union
  of the $g_iD(f)(k')$.  On the other hand, $G(k')x$ contains
  $D(f)(k')$ and hence also contains each $g_iD(f)(k')$; this proves
  the containment $\supseteq$ of \eqref{eq:orbit-claim} and completes
  the proof of the lemma.
\end{proof}

We are now going to show:

\begin{prop}
  \label{prop:finite-main}
  Let $k'$ be an algebraically closed extension field of $k$. Assume
  that $G$ has $n < \infty$ orbits on $W = W(k)$. Then each
  $G(k')$-orbit in $W(k')$ has a $k$-rational point. In particular,
  $G(k')$ has $n$  orbits on
  $W(k')$.  
\end{prop}

Note that if there are an infinite number of $G$-orbits on $W$, there
may indeed by $G(k')$-orbits on $G(k')$ without $k$ rational points.
This phenomenon already occurs in case $G$ acts trivially on a
positive dimensional variety $W$.

In view of Lemma \ref{lem:orbit-claim} and the fact that any $G$-orbit
in $W$ is a locally closed subvariety \cite{springer-LAG}*{Lemma 2.3.3}, it is
clear that Proposition \ref{prop:finite-main} follows from the Lemma
which follows.

\begin{lem}
  \label{lem:k'-union}
  Suppose that the irreducible affine $k$-variety $W$ is a union
  \begin{equation}
    \label{eq:union}
    W = L_1 \cup \cdots \cup L_n
  \end{equation}
  where the $L_i$ are non-empty, locally closed subvarieties, and that
  $k'$ is any field containing $k$. Then
  \begin{equation*}
    W(k') = L_1(k') \cup \cdots \cup L_n(k').
  \end{equation*}
\end{lem}

\begin{proof}  
  After possibly increasing $n$ and replacing the $L_i$ by smaller
  locally closed subvarieties, we may suppose for $i=1,2,\dots,n$ that
  the closure of $L_i$ is the closed set $\V(J_i)$ defined by an ideal
  $J_i = \sqrt{J_i} \ideal k[W]$, and that $L_i$ has the form $D(f_i)
  \cap \V(J_i)$ for a non-0 regular function $f_i \in k[W]$.
  
  The condition (\ref{eq:union}) may be restated: 
  \begin{itemize}
  \item[(*)] for each $k$-algebra homomorphism $\alpha:k[W] \to k$,
    there is some $1 \le i \le n$ with $\alpha(J_i) = 0$ and
    $\alpha(f_i) \not = 0$.
  \end{itemize}
  
  Any point of $W(k')$ is given by a $k$-homomorphism $\alpha:k[W] \to
  k'$.  To prove the lemma, we only must argue that $\alpha(J_i)=0$
  and $\alpha(f_i) \not = 0$ for some $1 \le i \le n$, since then
  $\alpha$ determines a point of $L_i(k')$.
  
  Let $I = \ker \alpha$. The algebra $k[W]/I$ is isomorphic to a
  $k$-subalgebra of the field $k'$; in particular, $I$ is a prime
  ideal and so the closed subset $\V(I)$ of $W$ is an
  irreducible $k$-variety.  
 
  Since $\V(I)$ is contained in the union of the $L_i$, it follows
  from Lemma \ref{lem:open} that $\W = \V(I) \cap L_m$ is a
  non-empty open subset of $\V(I)$ for some $1 \le m \le n$.  If $d =
  \dim \V(I)$, the closure of $\W=\V(I) \cap L_m$ is a closed subset
  of $\V(I)$ of dimension $d$; by irreducibility, $\V(I)$
  is precisely the closure of $\W$.  On the other hand, the closure of
  $\W$ lies in the closure of $L_m$, which is $\V(J_m)$; from this we
  find that $\V(I) \subseteq \V(J_m)$. Since $I = \sqrt{I}$ and $J_m
  = \sqrt{J_m}$ we deduce from Hilbert's Nullstellensatz that $J_m
  \subseteq I$; thus $\alpha(J_m) = 0$.  Since $\V(I) \cap L_m$ is
  non-empty, in particular $\V(I) \cap D(f_m)$ is non-empty; thus
  $f_m \not \in I$.  This means that $\alpha(f_m) \not = 0$, which
  completes the proof of the lemma.
\end{proof}

\begin{rem}
  A different proof of Proposition \ref{prop:finite-main} due to
  Guralnick may be found in \cite{gur}*{Prop. 1.1}.
\end{rem}

\subsection{Subvarieties of a linear $G$-representation}

Let $V$ be a finite dimensional $k$-vector space on which the
algebraic group $G$ acts linearly.  Let $W \subset V$ be an
irreducible $G$-invariant subvariety on which $G$ has finitely many
orbits. Assume as well that $kx \subset W$ for each $x \in W$.

Since the set $k^\times x$ lies in $W$, it only meets a finite number
of $G$-orbits; thus there is an orbit $\Orbit \subset W$ such that
$k^\times x \cap \Orbit$ is infinite.  Hence there is some $\beta \in
k^\times$ such that $k^\times x \cap G(\beta x)$ is infinite. Since
$G$ acts linearly on $V$, it follows at once that $k^\times x \cap Gx$
is infinite.

Consider the subgroup $N(x) = \{g \in G \mid gx \in k^\times x\} \le
G$; there is a homomorphism $\lambda:N(x) \to \Gm$ determined by the
condition $gx =\lambda(g)x$ for $g \in N(x)$.  Observe that the image
of $\lambda$ is an infinite subgroup of $\Gm$. Indeed, any $\alpha \in
k^\times$ such that $\alpha x \in k^\times x \cap Gx$ lies in the
image of $\lambda$.

Since $\Gm$ is a connected subgroup of dimension 1, the image of
$\lambda$ is in fact all of $\Gm$.  We conclude:
\begin{equation}
  \label{eq:Gm-orbit}
  \text{if $x \in W$, then $k^\times x \subset Gx$.}
\end{equation}

Fix $v,w \in W$, and assume that
\begin{equation*}
  av + bw \in W \ \text{for each $(a,b) \in k^2$}.
\end{equation*}
Since $W$ is stable under the scalar $k^\times$ action on $V$, this is
a ``projective'' condition; i.e. we may make instead the equivalent
assumption:
\begin{equation}
  \label{eq:vw-assumption}
   av + bw \in W \ \text{for each point}\  (a:b)   \in \Proj^1.
\end{equation}

\begin{prop}
  \label{prop:generic-point-orbit}
  Let $v,w\in W$ and assume that \eqref{eq:vw-assumption} holds.
  \begin{enumerate}
  \item There is a $G$-orbit $\Orbit \subset W$ and a non-empty open
    subset $\U \subset \Proj^1_{/k}$ such that $av + bw \in \Orbit$ if
    $(a:b) \in \U$ and $\dim G(av + bw) < \dim \Orbit$ if $(a:b) \in
    \Proj^1 \setminus \U$.
  \item Let $k' \supset k$ be an extension field and let $t\in k'$ be
    transcendental over $k$. Then $v + tw \in \Orbit(k')$ so that
    $\Orbit(k_1) = G(k_1)(v + tw)$ for any algebraically closed field
    $k_1$ containing $k'$.
  \end{enumerate}
\end{prop}

\begin{proof}
  Let $\phi:\Aff^2 \to W$ be the morphism $(a,b) \mapsto av +
  bw$. The image of $\phi$ is a closed and irreducible subvariety $S$
  of $W$.  Since $G$ has finitely many orbits on $W$, it follows from
  Lemma \ref{lem:open} that $S \cap \Orbit$ is open in $S$ for a unique
  $G$-orbit $\Orbit \subset W$. Moreover, since $S$ is closed, it is
  contained in the closure $\overline{\Orbit}$ of $\Orbit$.
  
  Thus $\U_1 = \phi^{-1}(\Orbit \cap S)$ is an open subset of
  $\Aff^2$ with the property that $av + bw \in \Orbit$ whenever
  $(a,b) \in \U_1$ and $av+bw \in \overline{\Orbit} \setminus \Orbit$
  whenever $(a,b) \in \Aff^2 \setminus \U_1$. 
  
  To complete the proof of (1), view $\Aff^2 \setminus 0$ as a
  $\Gm$-bundle $\pi:\Aff^2 \setminus 0 \to \Proj^1$.  Since $\pi$ is a
  flat morphism of finite type, it is open -- e.g. by \cite{liu}*{Exerc. 4.3.9} --
  so that $\U = \pi(\U_1)$ is the desired open subset of $\Proj^1$.
 
  For (2), let $\eta \in \Proj^1$ be the generic point. Identify
  $k(t)$ with $k(\Proj^1)$, and view $$\bar \eta = (1:t) \in
  \Proj^1(k')$$
  as a geometric point over $\eta$. Since $\eta$ is a
  point of $\U$ (in the sense of schemes), we have $\bar \eta \in
  \U(k')$.  Thus $v + tw \in \Orbit(k')$, and the remainder of (2)
  follows from Lemma \ref{lem:orbit-lemma}.
\end{proof}

\begin{rem}
  In the sequel, we will apply the previous result to $G=\GL(V)$
  acting by the adjoint representation on its Lie algebra $\gl(V)$.
  The nilpotent variety $\NN \subset \gl(V)$ satisifes $kX \subset
  \NN$ for each $X \in \NN$, and $\GL(V)$ has finitely many orbits on
  $\NN$. Moreover, \eqref{eq:vw-assumption} holds for any pair $X,Y
  \in \NN$ for which $[X,Y]=0$.
\end{rem}

\section{Background for $\GL(V)$}
\label{sec:background}

Let us recall how to recognize the unipotent radical of the
centralizer of a nilpotent element for the group $G=\GL(V)$.  If $A
\in \gl(V)$ is any nilpotent element, the $A$-exponent of $v \in V$ is
the non-negative integer 
\begin{equation}
  \label{eq:mu-defined}
  \mu(v)=\mu(A;v)=\min(n \ge 0 \mid A^nv=0).
\end{equation}
The vectors $v_1,\dots,v_n \in V$ are said to be $A$-independent
provided that the set
\begin{equation*}
  (*) \quad \{A^jv_i \mid 1 \le i \le n,\ 0 \le j \le \mu(v_i)-1\}
\end{equation*}
is linearly independent over $k$.  The vectors $v_1,\dots,v_n \in V$
form an $A$-basis if $(*)$ is a $k$-basis for $V$.

We recall some basic results. If $A \in \gl(V)$ is nilpotent, there is
an $A$-basis of $V$. If $v_1,\dots,v_n$ is an $A$-basis, ordered such
that $\mu(v_1) \ge \mu(v_2) \ge \cdots \ge \mu(v_n)$, write $\lambda =
(\lambda_1\ge \lambda_2 \ge \cdots \ge \lambda_n)$ for the partition
of $\dim V$ whose parts are $\lambda_i = \mu(v_i)$.  The
partition $\lambda$ is independent of the choice of $A$-basis for $V$,
and the $\GL(V)$-orbit of $A$ depends only on the partition $\lambda$,
which is thus called the partition of $A$.

A cocharacter of an algebraic group $G$ is a homomorphism $\Gm \to G$;
cocharacters of $\GL(V)$ may be identified with $\Z$-gradings of $V$.
Indeed, if $\chi:\Gm \to \GL(V)$ is a cocharacter, the weight spaces
\begin{equation*}
  V(m)=V(\chi;m) = \{v \in V \mid \chi(s)v = s^mv \ \forall s \in \Gm\}
\end{equation*}
determine a $\Z$-grading $V = \bigoplus_{m \in \Z} V(m)$ of $V$. Conversely, if
$V = \bigoplus_{m \in \Z} V(m)$ is a $\Z$-grading, there
is a unique cocharacter $\chi:\Gm \to G$ for which $V(m) = V(\chi;m)$.

We record:
\begin{lem}
  \label{lem:A-basis->cochar}
  Let $A \in \gl(V)$ be nilpotent with partition $\lambda$.  An
  $A$-basis $\{v_1,\dots,v_n\}$ of $V$ determines a unique cocharacter
  $\chi:\Gm \to \GL(V)$ for which $V(m)= V(\chi;m)$ is spanned by the
  vectors
  \begin{equation}
    \label{eq:grading-def}
    A^jv_i \quad \text{with} \quad m = -\lambda_i + 1 + 2j
  \end{equation}
  for $m \in \Z$.
\end{lem}

Note that the cocharacter $\chi$ depends only on $A$ and the choice of
an $A$-basis for $V$; we say that $\chi$ \emph{is a cocharacter associated to}
$A$.

Under the adjoint action of $\GL(V)$ on its Lie algebra $\gl(V)$, the
grading determined by $\chi$ has homogeneous components
\begin{align*}
  \gl(V)(m) &= \{C \in \gl(V) \mid C(V(j)) \subset V(j+m) \ \text{for each $j$}\} \\
  &= \{C \in \gl(V) \mid \Ad(\chi(s))C = s^mC  \ \forall s \in \Gm\}.
\end{align*}
for $m \in \Z.$ In particular, $A \in \gl(V)(2)$.

The cocharacter $\chi$ determines a unique parabolic subgroup
$P(\chi) < \GL(V)$ whose Lie algebra is
\begin{equation*}
  \plie(\chi) = \sum_{j \ge 0} \gl(V)(j);
\end{equation*}
moreover, if $U=R_uP(\chi)$, then
\begin{equation*}
  \ulie = \Lie(U) = \sum_{j > 0} \gl(V)(j).
\end{equation*}

\begin{prop}
  \label{centralizer-prop}
  Let $A \in \gl(V)$ be nilpotent. If $B \in \gl(V)$ satisfies
  $[A,B]=0$, then $B \in \plie(\chi)$, where the cocharacter $\chi$ is
  associated with $A$. Similarly, if $g \in \GL(V)$ satisfies
  $\Ad(g)A=A$, then $g\in P(\chi)$.
\end{prop}

\begin{proof}
  See \cite{jantzen-nil}*{3.10}.
\end{proof}

\begin{prop}
  \label{prop:cochar-conj}
  Any two cocharacters associated with $A$ are conjugate by an element
  of $\GL(V)$ centralizing $A$.
\end{prop}

\begin{proof}
  Indeed, any two $A$-bases are conjugate by an element centralizing
  $A$.
\end{proof}

\begin{prop}
  Let $A \in \gl(V)$ be nilpotent, and let $\chi$ be a cocharacter
  associated with $A$.  If $P=P(\chi)$, then $P$ is the instability
  parabolic subgroup for the unstable vector $A \in \glie$ in the
  sense of Kempf \cite{kempf-instab}. In particular, $P$ is
  independent of the choice of $A$-basis for $V$.
\end{prop}

\begin{proof}
  The fact that $P$ is the instability parabolic follows from the
  discussion (and references) in \S \ref{section:otherG}; see also
  \cites{jantzen-nil,premet,mcninch-rat}. The fact that $P$ is
  independent of the choice of $A$-basis for $V$ follows from general
  results about the instability parabolic. However, there is an
  elementary proof that $P$ is independent of the choice of $A$-basis:
  if $\chi$ and $\chi'$ are two cocharacters associated with $A$, then
  by Proposition \ref{prop:cochar-conj}, the cocharacters $\chi$ and
  $\chi'$ are conjugate by $g \in \GL(V)$ with $\Ad(g)A = A$. Thus
  $P(\chi') = gP(\chi)g^{-1} = P(\chi)$ since $g \in P(\chi)$ by
  Proposition \ref{centralizer-prop}.
\end{proof}

\begin{rem}
  If $L \supset K$ is a field extension and if $A \in \gl(V)(L)$ is
  nilpotent, then there is an $A$-basis $v_1,\dots,v_n \in V(L)$. For
  such a choice of $A$-basis, the homogeneous components $V(m)$ and
  $\gl(V)(m)$ are defined over $L$ for $m \in \Z$.  Equivalently: the
  cocharacter $\chi$ determined by this choice of $A$-basis is defined
  over $L$.  Thus the parabolic subgroup $P(\chi)$ is defined over $L$.
\end{rem}

The choice of cocharacter $\chi$ associated with $A$ determines a Levi
factor $L(\chi)$ in $P(\chi)$: take $L(\chi)$ to be the subgroup
$\prod_{i \in \Z} \GL(V(\chi;i)) \le \GL(V)$.

Denote by $C$ the centralizer of the nilpotent $A \in \GL(V)$, and
choose a cocharacter $\chi$ associated with $A$. We have:
\begin{prop}
  \label{prop:centralizer-desc}
  Let $C_\chi = C \cap L(\chi)$ and $R = C \cap R_uP(\chi)$. Then $C =
  C_\chi \cdot R$ is a semidirect product, $C_\chi$ is a reductive
  group isomorphic to a product of groups $\GL_r$ for various $r$,
  and $R$ is the unipotent radical of $C$.
\end{prop}

\begin{proof}
  \cite{jantzen-nil}*{Prop. 3.10 and Prop. 3.8.1}.
\end{proof}

\section{The main result}

We begin with a few preliminary results.

\subsection{Modifying an $A$-basis.}

Let $A$ be a nilpotent endomorphism of $V$, choose an $A$-basis
$\{v_1,\dots,v_n\}$ of $V$; put $\lambda_i = \mu(v_i)$ for $1 \le i
\le n$, and assume that $\lambda_1 \ge \cdots \ge \lambda_n$.

Assume that $B$ is a second nilpotent endomorphism of $V$ and that
$[A,B] = 0$.  The choice of $A$-basis made above determines a
cocharacter $\chi$ as in Lemma \ref{lem:A-basis->cochar}. By
Proposition \ref{centralizer-prop}, we may write $B = \sum_{i \ge 0}
B_i$ with $B_i \in \glie(V)(\chi;i)$.

Since $\chi(\Gm)$ normalizes the centralizer of $A$, we find that
$[A,B_0] = 0$ as well. It follows that the endomorphism $B_0$ is
determined by its values on the $A$-basis vectors $v_1,\dots,v_n$.  In
particular, if $B_0 \ne 0$, then $B_0v_i \ne 0$ for some $1 \le i \le
n$.

 \begin{lem}
   \label{lem:change-basis}
   Fix $1 \le i \le n$, and assume that $B_0v_i \ne 0$. Then 
   \begin{enumerate}
   \item[(1)] $\mu(Bv_i) = \mu(A;Bv_i) = \lambda_i$, and
   \item[(2)] for some $j \ne i$ with $1 \le j \le n$ and $\lambda_j =
     \lambda_i$, the vectors
   \begin{equation*}
       v_1,\dots,v_{j-1},Bv_i,v_{j+1},\dots,v_n
     \end{equation*}
     form an $A$-basis for $V$.
   \end{enumerate}
 \end{lem}

 \begin{proof}
   Since $A$ and $B$ commute, it is clear that $A^{\lambda_i}Bv_i =
   0$. To complete the proof of (1), we must argue that $A^{\lambda_i
     -1}Bv_i \ne 0$.  According to \cite{jantzen-nil}*{3.1(1)}, we
   have
   \begin{equation}
          \label{eq:B-expression}
          Bv_i = \sum_{j=1}^n \sum_{\ell = 
            \max(0,\lambda_j - \lambda_i)}^{\lambda_j - 1} c_{\ell,j} A^\ell v_j
   \end{equation}
   for certain $c_{\ell,j} \in k$.
   It follows that
   \begin{equation}
     \label{eq:B_0-expression}
     B_0v_i = \sum_{\lambda_j = \lambda_i} a_j v_j \quad \text{with}\ a_j = c_{0,j} \in k;
   \end{equation}
   moreover, with notation as in
   \eqref{eq:B-expression}
   \begin{equation}
     \label{eq:B-expression-refined}
     \begin{split}
       & Bv_i = B_0v_i + w + Ax \\ & \text{where $w = \sum_{\lambda_j <
           \lambda_i} c_{0,j}v_j$ so that $A^{\lambda_i - 1}w = 0$ and
         $A^{\lambda_i +1}x = 0$.}
     \end{split} 
   \end{equation} 
   Indeed, to verify \eqref{eq:B-expression-refined}, notice that if
   $\lambda_j < \lambda_i$, then $A^{\lambda_i-1}v_j = 0$ so that
   $A^{\lambda_i-1}w = 0$. Now notice that $Bv_i - B_0v_i - w$ has the
   form $Ax$ for some $x \in V$.  Finally, since $Bv_i$, $B_0v_i$ and
   $w$ lie in the kernel of $A^{\lambda_i}$, so does $Ax$.

   It follows that
   \begin{equation*}
     A^{\lambda_i-1}Bv_i \equiv A^{\lambda_i-1}B_0v_i \pmod{A^{\lambda_i}V}.
   \end{equation*}
   Since $B_0v_i$ is non-zero and is a linear combination of the $v_k$
   with $\lambda_k = \lambda_i$, it is clear that $A^{\lambda_i -
     1}B_0v_i \not \equiv 0 \pmod{A^{\lambda_i} V}$; thus
   $A^{\lambda_i - 1}Bv_i$ is non-zero $\pmod{A^{\lambda_i} V}$. In
   particular, $A^{\lambda_i - 1}Bv_i$ is non-zero, which completes
   the proof of (1).

   As to (2), one knows that $B_0$ is nilpotent since $B$ is
   nilpotent. Thus the vectors $v_i$ and $B_0v_i$ are linearly
   independent. In the above expression \eqref{eq:B_0-expression} for
   $B_0v_i$, it follows that $a_j \ne 0$ for some $1 \le j \le n$ with
   $j \ne i$ and $\lambda_j = \lambda_i$.

   We are going to prove that (2) holds for this value of $j$. As a
   preliminary step, notice that $\{v_1,\dots,v_n\}$ remains an
   $A$-basis if we replace $v_j$ by $B_0 v_i$; thus we may and will
   suppose that $B_0v_i = v_j$.

   Let us write $\lambda = \lambda_i = \lambda_j$.  With notation as
   in \eqref{eq:B-expression-refined}, recall that
   $A^{\lambda-1} w = 0.$
   Let now
   \begin{equation*}
     u_s = \left\{
       \begin{matrix}
         v_s & s \ne j \\
         Bv_i & s=j.
       \end{matrix}
       \right .
     \end{equation*}
     We must show that $\{u_1,\dots,u_n\}$ is an $A$-basis of $V$. To
     see this, let $f_1,f_2,\dots,f_n \in K[z]$ be polynomials for
     which $\displaystyle \sum_{s=1}^n f_s(A) u_s = 0.$ We must argue
     that $f_s$ is divisible by $z^{\lambda_s}$ for each $1 \le s \le
     n$.  In fact, it is enough to argue that $z^{\lambda}$ divides
     $f_j$, since then the result follows from the $A$-independence of
     the set $\{v_1,\dots,v_{j-1},v_{j+1},\dots,v_n\}$.

   Using \eqref{eq:B-expression-refined}
   we have
  \begin{equation*}
     0 = f_j(A)Bv_i + \sum_{s \ne j} f_s(A)v_s  =  
     f_j(A)v_j + f_j(A)w + f_j(A)Ax + \sum_{s \ne j} f_s(A)v_s.
   \end{equation*}
   If $f_j = 0$, then of course $z^{\lambda}$ divides $f_j$ and the
   proof is complete.  If $f_j \ne 0$, let $\mu \ge 0$ be maximal such
   that $z^\mu \mid f_j$, and write $f_j = z^\mu \cdot g$ for a
   polynomial $g \in K[z]$ having non-zero constant term.  We find
   then that
   \begin{equation*}
     A^\mu g(A)v_j \equiv -A^\mu g(A) w - \sum_{s \ne j} 
     f_s(A)v_s \pmod{A^{\mu+1}V}.
   \end{equation*}
   Since $w$ is a linear combination of $v_s$ with $\lambda_s <
   \lambda_j$, the right hand side is congruent to an expression of
   the form $\displaystyle \sum_{s \ne j} h_s(A) v_s$ modulo
   $A^{\mu+1}V$ for polynomials $h_s \in K[z]$.  Since the vectors
   $v_1,\dots,v_n$ are $A$-independent, it follows that $A^\mu g(A)v_j
   \equiv 0 \pmod {A^{\mu+1}V}$.  Since $g$ has non-zero constant
   term, this is only possible if $\mu \ge \lambda$, as required.
 \end{proof}

 \subsection{Recognizing the partition of a nilpotent endomorphism
   $A$}

  Let $A$ be a nilpotent endomorphism of $V$. Let $v_1,\dots,v_n \in
  V$ and let $\lambda = (\lambda_1 \ge \lambda_2 \ge \cdots \ge
  \lambda_n)$ be a partition of $\dim V$. Suppose that the set of
  vectors
  \begin{equation*}
    \B=\{ A^j v_i \mid 1 \le i \le n, \ 0 \le j \le \lambda_i  -1\}
  \end{equation*}
  forms a $k$-basis for $V$. For $i \ge 0$ write $V_\ell$ for the
   span of $\{A^jv_i \mid 1 \le i \le \ell,\
  0 \le j < \lambda_i\}$; thus $V_0 = 0$.  

\begin{lem}
  \label{lem:partition-condition}
  The following are equivalent:
  \begin{enumerate}
  \item   $A^{\lambda_j}v_j \in A^{\lambda_j}V_{j-1}$ for each $1 \le j \le n$,
  \item for $1 \le j \le n$ there are vectors $w_j \in V_{j-1}$ such
    that $A^{\lambda_j}(v_j - w_j) = 0$ for $1 \le j \le n$ and such
    that $\{v_j - w_j \mid 1 \le j \le n\}$ is an $A$-basis of $V$,
  \item  $\lambda$ is the partition of $A$.
  \end{enumerate}
  In particular, if $\lambda$ is the partition of $A$, then each
  subspace $V_\ell$, $1 \le \ell \le n$, is $A$-invariant.
\end{lem}

\begin{proof}
  To prove $(1) \implies (2)$, choose for each $j=1,2,\dots,n$ a vector
  $w_j \in V_{j-1}$ for which
  \begin{equation*}
    A^{\lambda_j}v_j = A^{\lambda_j}w_j \quad \text{hence} \quad
    A^{\lambda_j}(v_j - w_j) = 0.   
  \end{equation*}
  To see that the vectors $v'_j=v_j - w_j$ for $1 \le j \le n$ form an
  $A$-basis, just note that if $M$ is the matrix of coefficients
  obtained upon expressing the vectors $A^sv_t'$ in terms of the
  $K$-basis $\{A^lv_m\}$, then $M$ is unipotent and hence invertible.

  The assertions $(2) \implies (1)$ and $(2) \implies (3)$ are immediate.
  
  We finally prove $(3) \implies (2)$. Since $\lambda$ is the
  partition of $A$, $A^{\lambda_1} = 0$; in particular,
  $A^{\lambda_1}v_1 = 0$. Apply \cite{lang}*{Lemma III.7.6} to see
  that $(\lambda_2 \ge \cdots \ge \lambda_n)$ is the partition of the
  nilpotent endomorphism $\overline{A}$ of $V/V_1$ induced by $A$; by
  induction on $n$ we find vectors $w_j' \in V_{j-1}$ for $2 \le j \le
  n$ such that
  \begin{equation*}
    A^{\lambda_j}(v_j - w_j') \in V_1 \quad \text{for}\ 2 \le j \le n
  \end{equation*}
  and such that $v_2-w_2',\dots,v_n-w_n'$ is an $\overline{A}$-basis
  of $V/V_1$.  Another application of \cite{lang}*{Lemma III.7.6} now
  gives vectors $w''_2,\dots,w''_n \in V_1$ for which 
  \begin{equation*}
    v_1,v_2- w_2'-w_2'',\dots,v_n-w_n'-w_n''
  \end{equation*}
  is an $A$-basis for
  $V$. Since $V_1 \subset V_{j-1}$ for $j \ge 2$, we have $w_j=w_j'-w_j''
  \in V_{j-1}$ as desired; thus (2) indeed holds.
\end{proof}

\subsection{A nilpotent element of $\gl(V)$ over $K(t)$}

Let $p$ denote the characteristic of $K$, and recall that $t$ is
transcendental over $K$.  Let us fix nilpotent elements $X,Y \in
\gl(V)(K)$, and let us suppose that $[X,Y]=0$.

Write $\lambda = (\lambda_1 \ge \lambda_2 \ge \cdots \ge \lambda_n)$
for the partition of $X$, and fix once and for all an $X$-basis
$v_1,\dots,v_n \in V(K)$ for $V$. Thus
\begin{equation*}
  \B^0=\{X^j v_i \mid 1 \le i \le n,\ 0 \le j \le \lambda_i -1 \}
\end{equation*}
is a $K$-basis for $V(K)$.
  
Consider the localization $\A=K[t]_{(t)}$ of the polynomial ring
$K[t]$ at the maximal ideal $tK[t]$; its field of fractions is
$F=K(t)$, and its maximal ideal is $\mm = (t)= t\A$. Write $\VV = V(K)
\tensor_K \A$. Each of the vectors in the set
\begin{equation*}
  \B^t=\{ (X+tY)^j v_i \mid 1 \le i \le n, \ 0 \le j \le \lambda_i -1\}
\end{equation*}
lies in $\VV$. By assumption, the image in $V(K) = \VV/t\VV$ of
$\B^t$ is $\B^0$; by the Nakayama lemma, $\B^t$ forms an $\A$-basis
for $\VV$. In particular, $\B^t$ is an $F=K(t)$-basis for $V(F)$.
  
For each $1 \le \ell \le n$, let us write $V^0_\ell(K)$ for the
$K$-subspace of $V(K)$ spanned by the vectors
\begin{equation*}
  \B^0_\ell = \{X^jv_i \mid 1 \le i \le \ell,\ 0 \le j \le \lambda_i -1\}.
\end{equation*}
Similarly, let $V^t_\ell(F)$ be the $F$-subspace of $V(F)$ spanned by the vectors
\begin{equation*}
  \B^t_\ell = \{(X+tY)^jv_i \mid 1 \le i \le \ell,\ 0 \le j \le \lambda_i -1\},
\end{equation*}
and let $\VV_\ell$ be the $\A$-submodule of $\VV$ spanned by
$\B^t_\ell$. Of course, the image of $\VV_\ell$ in $V(K) = \VV/t\VV$
is $V^0_\ell$.

\begin{lem}
  \label{lem:direct-summand}
  For $1 \le \ell \le n$, $\VV_\ell$ is a direct summand of $\VV$ as
  an $\A$-module. We have in particular:
  \begin{enumerate}
    \item $\VV_\ell = V^t_\ell(F) \cap \VV$, and
    \item  $t\VV_\ell = \VV_\ell \cap t\VV$.
  \end{enumerate}
\end{lem}

\begin{proof}
  Since $\B^t$ is an $\A$-basis of $\VV$, the lemma is immediate.
\end{proof}

\begin{lem}
  \label{lem:tech}
  Assume that the partition of $X+tY$ coincides with that of $X$; i.e.
  assume that $X + tY$ and $X$ are $\GL(V)$-conjugate. For each $1 \le
  \ell \le n$, we have:
  \begin{enumerate}
  \item $\VV_\ell$ is $X+tY$-invariant, and
  \item $(X+tY)^{\lambda_\ell} v_\ell \in \VV_{\ell-1}.$
  \end{enumerate}
\end{lem}

\begin{proof}
  Fix $1 \le \ell \le n$. Since $\lambda$ is the partition of $X+tY$,
  Lemma \ref{lem:partition-condition} shows that each $V^t_\ell(F)$ is
  $X+tY$-invariant.  Since $\VV_\ell = V^t_\ell(F) \cap \VV$ by Lemma
  \ref{lem:direct-summand}(1), the $X+tY$-invariance of $\VV_\ell$
  results from that of $\VV$ and of $V^t_\ell(F)$; this proves (1).
  
  Since
  \begin{equation*}
    (X+tY)^{\lambda_\ell}v_\ell \in  (X+tY)^{\lambda_\ell}V^t_{\ell-1}(F) 
    \subset V^t_{\ell-1}(F),
  \end{equation*}
  we have $(X+tY)^{\lambda_\ell}v_\ell \in \VV \cap V^t_{\ell-1}(F)$
  by another application of Lemma \ref{lem:direct-summand}(1).
\end{proof}

\begin{prop}
  \label{prop:Y0-vanish}
  Assume that the partition of $X+tY$ coincides with that of $X$; i.e.
  assume that $X + tY$ and $X$ are $\GL(V)$-conjugate. Let $\chi$ be
  the $K$-cocharacter associated with $X$ determined by the $X$-basis
  $v_1,\dots,v_n$, and write $Y = Y_0 + Y_+$ with
  \begin{equation*}
    Y_0 \in \gl(V)(\chi;0) \quad \text{and} \quad  
    Y_+ \in \sum_{j>0} \gl(V)(\chi;j).
  \end{equation*}
  If $p>0$ assume that $X^{p-1} = 0$. Then $Y_0=0$.
\end{prop}

\begin{proof}
  We assume that $Y_0 \ne 0$ and deduce a contradiction.  Let $1 \le
  \ell \le n$ be minimal with $Y_0v_\ell \ne 0$.  After possibly
  re-ordering those members of the $X$-basis $v_1,v_2,\dots,v_n$ for
  which $\lambda_k = \lambda_\ell$, we may suppose that $\lambda_k >
  \lambda_\ell$ whenever $k < \ell$. According to Lemma
  \ref{lem:change-basis}, we may and will assume that $Yv_\ell = v_j$
  for some $j > \ell$ with $\lambda_j = \lambda_\ell$.
  
  Since $\lambda$ is the partition of $X+tY$, Lemma \ref{lem:tech}
  shows that $(X+tY)^{\lambda_\ell}v_\ell \in \VV_{\ell-1}$. Since
  $X^{\lambda_\ell}v_\ell = 0$, we find by Lemma \ref{lem:direct-summand}(2):
  \begin{equation*}
    (X+tY)^{\lambda_\ell}v_\ell \in \VV_{\ell-1} \cap t\VV = t\VV_{\ell -1}.
  \end{equation*}
  Thus we see 
  \begin{equation*}
    \dfrac{1}{t}(X+tY)^{\lambda_\ell}v_\ell
    = \sum_{j=1}^{\lambda_\ell} t^{j-1}\dbinom{\lambda_\ell}{j}
    X^{\lambda_\ell - j}Y^j v_\ell \in \VV_{\ell-1}.
  \end{equation*}
  Since the image of $\VV_{\ell-1}$ in the quotient $V(K) = \VV/t\VV$
  is $V^0_{\ell-1}(K)$, it follows that
  \begin{equation*}
    \lambda_\ell X^{\lambda_\ell-1}Yv_\ell = \lambda_\ell X^{\lambda_\ell-1}v_j 
    \in V^0_{\ell-1}(K).
  \end{equation*}
  If $p>0$, the condition $X^{p-1} = 0$ shows that $\lambda_\ell < p$;
  so in every case, $\lambda_\ell$ is non-zero in $K$. It follows that
  $X^{\lambda_\ell -1}v_j=X^{\lambda_j - 1} v_j \in V^0_{\ell-1}(K)$,
  contradicting the assumption that $v_1,\dots,v_n$ is an $X$-basis
  for $V$. This completes the proof.
\end{proof}

\subsection{A nilpotent element of $\gl(V)$ over $\Proj^1$}

Let $X,Y \in \gl(V)(K)$ be nilpotent with $[X,Y]=0$, and let $\OO$
denote the structure sheaf of $\Proj^1=\Proj^1_{/K}$. Write $\LL =
V(K) \tensor_K \OO$, so that $\LL$ is a free sheaf of $\OO$-modules on
$\Proj^1$. If $\eta$ denotes the generic point of $\Proj^1$, the stalk
$\OO_\eta = K(\Proj^1)$ identifies with $F=K(t)$, and the stalk
$\LL_\eta$ identifies with $V(F)$. 

Choose an $A=X+tY$-basis $v_1,\dots,v_n \in V(F)$; for $1 \le i \le n$
and $j \ge 0$, we may regard $A^iv_j$ as an element of $\LL_\eta$.
Thus we may choose an affine open subset $\W \subset \Proj^1$ such
that $t$ is regular on $\W$ and such that $A^jv_i \in \Gamma(\W,\LL)$
for $1 \le i \le n$ and $0 \le j < \lambda_j$.

For a point $x \in \Proj^1$, denote by $\mm_x$ the maximal ideal of
the stalk $\OO_x$, and let $K(x)$ be the field of fractions of
$\OO_x/\mm_x$; the $K(x)$-vector space $\LL_x \tensor_{\OO_x} K(x)$
may be identified with $V(K(x)) = V(K) \tensor_K K(x)$.  If $\bar{x} =
(a:b)$ is a geometric point of $\W$ over $x$, then $X,Y,X+tY$ act on
$\LL_x$ and so on $V(K(x))$; the maps induced on $V(K(x))$ are
respectively $X,Y$, and some non-zero multiple of $aX + bY$
\footnote{The geometric point $\bar{x} = (a:b)$ over $x$ is determined
  by a field embeding $\iota:K(x) \to L$ for a separably closed field
  $L$. We have assumed that $t$ is regular at $x$ - i.e. $a \ne 0$ so that
  $t \in \OO_x$; if $\bar{t}$ denotes the image in $K(x)$ of $t \in
  \OO_x$, then $\iota(\bar{t})$ is a multiple of $b/a$. Now, $\iota$
  determines an embedding $\iota:V(K(x)) \to V(L)$; the map $aX + bY$
  leaves stable the image of $\iota$, and coincides with some multiple
  of $X+\bar{t}Y:V(K(x)) \to V(K(x))$. In this sense, $aX + bY$ is
  independent of the choice of geometric point $\bar{x}$.}.

We now have:

\begin{lem}
  \label{lem:open-A-basis}
  If $v_1,\dots,v_n \in V(F)$ is an $(X+tY)$-basis for $V$, there is a
  non-empty open subset $\U$ of $\Proj^1$ such that
  \begin{enumerate}
  \item $v_1,\dots,v_n \in \LL(\U)$, 
  \item the vectors $A^jv_i$ for $1 \le i \le n$ and $0 \le j <
    \lambda_i$ form a basis for $\LL(\U)$ over $\OO(\U)$, and
  \item for each $x \in \U$, the vectors $v_1,\dots,v_n \in V(K(x))$
    form an $(aX + bY)$-basis of $V$ for any geometric point $(a:b)$
    over $x$.
  \end{enumerate}
\end{lem}

\begin{proof}
  With notation as before, let $\MM = \bigwedge^{\dim V}\LL$, and
  consider the element
  \begin{equation*}
    \omega = \bigwedge_{j=1}^n \bigwedge_{i=0}^{\lambda_j -1} 
    A^iv_j \in \Gamma(\W,\MM).
  \end{equation*}
  Let $\U$ be a non-empty affine open subset of $\W$ for which the
  germ $\omega_x$ does not lie in $\mm_x\MM_x$ for all points $x \in
  \U$ [of course, the set of all $x \in \W$ having that property is
  non-empty and open].
  
  By construction, the vectors $\{A^jv_i \mid 1 \le i \le n, 0 \le j <
  \lambda_i\}$ form an $\OO(\U)$-basis of $\LL(\U)$, and the lemma
  follows.
\end{proof}

\subsection{The main theorem}

\begin{theorem}
  \label{theorem:GL(V)}
  Consider the nilpotent element $A=X + tY \in \gl(V)(F)$ where
  $X,Y \in \gl(V)(K)$ are nilpotent and $[X,Y]=0$. If $p>0$, assume
  that $A^{p-1} = 0$.  
  \begin{enumerate}
  \item  $X,Y \in \Lie R_uC$, where $C$ is the
    centralizer of $A = X+tY \in \gl(V)(F)$ in $\GL(V)$.
  \item There is a non-empty open subset $\U$ of $\Proj^1$ such that
    $X,Y \in \Lie R_uC_{(a:b)}$ for each geometric point $(a:b)$ of
    $\U$, where $C_{(a:b)}$ is the centralizer of $aX+bY$ in $\GL(V)$.
  \end{enumerate}
\end{theorem}

Before giving the proof, let me first give an example to demonstrate
that the theorem is not correct without some hypothesis on $A$.
\begin{example}
  Let $X' \in \gl(V)(K)$ be a regular nilpotent element, and write
  $d = \dim V$. Choose $v \in V(K)$ for which $\{v_i = (X')^iv \mid i
  = 0,1,\dots,d-1\}$ is a basis for $V$; we will write $v_i =
  (X')^iv$ for $i \ge 0$ so that $v_i=0$ for $i \ge d$.  Now let
  \begin{equation*}
    X = X' \oplus X' \in \gl(V \oplus V)
  \end{equation*}
  and let
  \begin{equation*}
    Y = ((v,w) \mapsto (0,v)) \in \gl(V \oplus V).
  \end{equation*}
  Of course, $[X,Y] = 0$. We set $A=X+tY \in \gl(V \oplus V)(F)$ 
  and write $C \le \GL(V \oplus V)$ for the centralizer of $A$.
  
  For $m \ge 0$,  we have
  \begin{equation*}
    A^m=(X + tY)^m = X^m + mtX^{m-1}Y.
  \end{equation*}
  If $w_1 = (v,0)$ and $w_2 = (0,v)$ we have:
  \begin{equation*}
    A^mw_1 = (v_m,mtv_{m-1}) \quad \text{and} \quad  A^mw_2 = (0,v_m)
  \end{equation*}
  for $m \ge 0$, where we have put $v_{-1}=0$.
  
  If $d \ne 0$ in $K$, the reader may verify that the partition of $A$
  is $\lambda = (d+1,d-1)$.  Since this partition has distinct parts,
  a Levi factor of $C$ is a torus so indeed
  $X,Y \in \Lie R_uP$.
  
  However, if $d = 0$ in $K$, then $A$ has partition $(d,d)$.  To see
  this, observe that $A^dw_1 = (0,dtv_{d-1}) = 0$, and $A^dw_2=0$; now
  verify that $w_1,w_2$ is an $A$-basis of $V\oplus V$. It is not true
  that $X \in \Lie R_uC$.  Indeed,
  \begin{equation*}
    Xw_1 = (v_1,0) = (v_1,tv_0) - t(0,v_0) = Aw_1 - tw_2;
  \end{equation*}
  since $w_1,w_2 \in V(-d+1)$, we find that $X_0 \not = 0$ so that $X
  \not \in \Lie R_uC$ [and so, of course, also $Y \not \in \Lie R_uC$].
  If $d=p$, $A$ is  $p$-nilpotent, i.e. we have $A^p =
  0$, but $A^{p-1} \not = 0$.
  \begin{center}
    * \qquad * \qquad  *
  \end{center}
\end{example}

\begin{proof}[Proof of Theorem  \ref{theorem:GL(V)}]
  First use Lemma \ref{lem:open-A-basis} to find an $(X+tY)$-basis
  $v_1,\dots,v_n$ for $V(F)$ and an open subset $\U \subset \Proj^1$
  satisfying the conclusion of that Lemma. If $x \in \U$ and $(a:b)$
  is a geometric point over $x$, the $aX+bY$-basis of $\LL_x
  \tensor_{\OO_x} K(x) = V(K(x))$ obtained from the $v_i$ determines a
  cocharacter $\chi_{(a:b)}$ associated to $aX+bY$.  Especially,
  $\chi_{(1:t)}$ is the cocharacter associated with $X+tY$ determined
  by the vectors $v_i \in V(F)$.
  
  Now, write $Y = Y_0 + Y_+$ for unique elements 
  \begin{equation*}
    Y_0 \in
  \gl(V)(\chi_{(1:t)};0) \quad\text{and}\quad Y_+ \in \sum_{j > 0}
  \gl(V)(\chi_{(1:t)};j).
  \end{equation*}
  Since the $\OO(\U)$-basis of $\LL(\U)$ determined by the $v_i$
  consists of weight vectors for the torus $\chi_{(1:t)}(\Gm)$, and
  since $Y(\LL(\U)) \subset Y(\LL(\U))$, one has that 
  \begin{equation*}
    Y_0(\LL(\U))
    \subset \LL(\U) \quad \text{and} \quad Y_+(\LL(\U)) \subset \LL(\U),
  \end{equation*}
  or -- what is the same -- one has that
  \begin{equation*}
    Y_0,Y_+  \in \gl(V)(\U) = \gl(V)(K) \tensor_K \OO(\U).
  \end{equation*}
  For each geometric point $(a:b)$ over $x \in \U$, write
  $(Y_0)_{(a:b)}$ and $(Y_+)_{(a:b)}$ for the images of $Y_0, Y_+$ in
  $\gl(V)(K(x)) = \gl(V)(\OO_x) \tensor_{\OO_x} K(x)$.  We have:
  \begin{equation*}
    (Y_0)_{(a:b)} \in \gl(V)(\chi_{(a:b)};0)
  \end{equation*}
  and
  \begin{equation*}
    (Y_+)_{(a:b)} \in \sum_{j>0} \gl(V)(\chi_{(a:b)};j).
  \end{equation*}
  Thus the theorem will follow from Proposition
  \ref{prop:centralizer-desc} provided that we only show $Y_0=0$.
  Moreover, it is enough to show that $(Y_0)_{(a:b)}=0$ for all
  geometric points $(a:b)$ in some dense subset of $\U$. 
  
  Writing $K(\Proj^1)= K(t)$, we may apply Proposition
  \ref{prop:generic-point-orbit} to find a non-empty open subset $\U'
  \subset \U$ such that $aX + bY$ is $\GL(V)$-conjugate to $X + tY$
  for each geometric point $(a:b)$ of $\U'$.
  
  We are now going to show that $(Y_0)_{(1:s)}=0$ for each point of
  $\U'$ of the form $(1:s)$ with $s \in K$. Since we may evidently
  replace $K$ by an algebraic extension, we may and will suppose that
  $K$ is infinite; thus such points are indeed dense in $\U'$ and
  hence in $\U$.
  
  So fix such a point $(1:s)$. Since $(1:s)$ is a point of $\U'$, we
  know that $X+sY \in \GL(V)$ is conjugate to $X+tY$. Since $t$ and
  $t+s$ are both transcendental over $K$, $X+sY$ and $X + (t+s)Y$ have
  evidently the same partition; thus $X+sY$ is conjugate to $X+(t+s)Y$
  as well. We may now apply Proposition \ref{prop:Y0-vanish} to the
  elements $X+sY$ and $Y$ to see that $(Y_0)_{(1:s)}=0$ as desired.
  This completes the proof of the theorem.
\end{proof}

\section{Other semisimple groups}
\label{section:otherG}

Consider now more general groups $G$; for ease of exposition I'll
assume that $G$ is semisimple over $K$, and that the characteristic of
$K$ is \emph{very good} for $G$.

Let $X \in \glie$ be nilpotent. A cocharacter $\phi:\Gm
\to G$ is \emph{associated to $X$} provided that:
\begin{itemize}
\item[A1.] $X \in \glie(\phi;2)=$ the 2-weight space of the torus
  $\phi(\Gm)$ under the adjoint representation on $\glie$, and
\item[A2.] for some choice of maximal torus $S < C_G(X)$, the image of
  $\phi$ lies in $(L,L)$, where $L$ is the Levi subgroup of $G$
  defined by $L=C_G(S)$.
\end{itemize}

When $G=\GL(V)$, the reader may easily check that the above
definition agrees with that given in \S \ref{sec:background}; namely,
if $X \in \gl(V)$ is nilpotent, then the cocharacters determined by
$X$-bases of $V$ as in Lemma \ref{lem:A-basis->cochar} are precisely
those satisfying A1 and A2.  For any $G$, the nilpotent element $X$ is
\emph{distinguished} in $\Lie(L)$ for a Levi subgroup $L$ as in A2;
for more on this see \cite{jantzen-nil}*{\S4--5}.

\begin{rem}
  \label{rem:JM}
  When $p=0$, the map $\tau \mapsto d\tau(1)$ is a bijection between
  cocharacters associated with $X$ and the set of all $H \in
  [X,\glie]$ such that $[H,X] = 2X$; cf. \cite{jantzen-nil}*{\S5.5}.
  Thus the cocharacters associated with $X$ are precisely those
  obtained by the Jacobson-Morozov Lemma.
\end{rem}

Under our assumptions on $G$, there are always cocharacters associated
to $X$ \cite{mcninch-rat}*{Prop. 16}; see also \cite{premet}.  If $X$
is $K$-rational, one can even find a cocharacter associated to $X$
which is defined over $K$; see \cite{mcninch-rat}*{Theorem 26}.  Any
cocharacter $\chi:\Gm \to G$ determines a parabolic subgroup $P(\chi)$
of $G$; namely, the unique parabolic whose Lie algebra is
$\bigoplus_{i \ge 0} \glie(\chi;i)$ where
\begin{equation*}
  \glie(\chi;i) = \{X \in \glie \mid \Ad(\chi(s))X = s^i X\ 
  \forall s \in k\}.
\end{equation*}
According to \cite{jantzen-nil}*{5.9}, the parabolic subgroup
$P(\phi)$ is independent of the choice of cocharacter $\phi$
associated to $X$; it is the instability parabolic of Kempf and
Rousseau \cite{mcninch-rat}*{Prop.  18}.

The analogues of Propositions \ref{centralizer-prop} and
\ref{prop:centralizer-desc} hold. Namely,
\begin{prop}
  \label{G-centralizer-prop}
  Let $A \in \glie(K)$ be nilpotent, let $\chi$ be a cocharacter
  associated with $A$, let $P=P(\chi)$, let $\plie = \Lie(P)$, and let $C=C_G(A)$ be the
  centralizer of $A$. Then
  \begin{enumerate}
  \item $C_G(A)$ is defined over $K$ and $\Lie C_G(A) = \clie_\glie(A)$,
  \item $\clie_\glie(A) \subset \plie$ and $C_G(A) \subset
    P$, and
  \item if $L(\chi)$ denotes the centralizer in $G$ of $\chi(\Gm)$,
    then $C_\chi = C \cap L(\chi)$ is reductive, $R_uC = C \cap R_uP$,
    and $C = C_\chi \cdot R_uC$ is a Levi decomposition.
  \end{enumerate}
\end{prop}

\begin{proof}
  (1) follows from the separability of orbits for semisimple groups in
  very good characteristic; see \cite{springer-steinberg}*{I.5.2 and
    I.5.6} together with \cite{springer-LAG}*{Prop. 12.1.2}.  (2) is
  \cite{jantzen-nil}*{Prop. 5.9}. (3) is \cite{jantzen-nil}*{Prop.
    5.10 and 5.11}; see also \cite{mcninch-rat}*{Corollary 29}.
\end{proof}

We want to consider the following hypothesis on $G$:
\begin{enumerate}
\item[\HypL] There is a representation $\rho:G \to \GL(V)$ defined
  over $K$ such that $d\rho$ is injective, and such that for each
  nilpotent $X \in \glie$ and each cocharacter $\chi$ of $G$
  associated with $X$, the cocharacter $\rho \circ \chi$ is of
  $\GL(V)$ associated with $d\rho(X) \in \gl(V)$.
\end{enumerate}

\begin{rem}
  It follows from Remark \ref{rem:JM} that the condition \HypL holds
  for any faithful representation $(\rho,V)$ when $\operatorname{char}
  K = 0$. Indeed, let $X \in \glie$ be nilpotent, let $\chi$ be a
  cocharacter of $G$ associated with $X$, and let $H = d\chi(1)$.
  Then $d\rho(H) = d(\rho\circ \chi)(1)$ in $\gl(V)$. Moreover,
  clearly $[d\rho(H),d\rho(X)] = d\rho([H,X]) = d\rho(2X) = 2d\rho(X)$
  and $d\rho(H) \in [d\rho(X),\gl(\glie)]$ so that $\rho \circ \chi$ is
  associated with $d\rho(X)$ by Remark \ref{rem:JM}.  This verifies
  \HypL.
\end{rem}

The following is an immediate consequence of Theorem
\ref{theorem:GL(V)}:

\begin{theorem}
  \label{theorem:general-G}
  Let $G$ be semisimple algebraic group defined over $K$, assume that
  the characteristic of $K$ is very good for $G$, and assume that
  \HypL holds. Let $X,Y \in \glie(K)$ with $[X,Y]=0$, and suppose that
  $d\rho(X+tY)^{p-1}=0$.
  \begin{enumerate}
  \item Then $X,Y \in \Lie R_uC$ where $C =
    C_G(X+tY)$ is the centralizer of $X+tY$.
  \item There is a non-empty open subset $\U$ of $\Proj^1$ such that
    for each geometric point $(a:b)$ of $\U$, we have $X,Y
    \in \Lie R_uC_{(a:b)}$, where $C_{(a:b)}=C_G(aX + bY)$.
  \end{enumerate}
\end{theorem}

\begin{lem}
  \label{lem:optimal}
  Let $X \in \glie$ satisfy $X^{[p]} = 0$, and suppose that $\chi$ is a
  cocharacter associated with $X$.
  \begin{enumerate}
  \item There is a homomorphism $\psi:\SL_2 \to G$ such that
      \begin{equation*}
        d\psi \left (
          \begin{matrix}
            0 & 1 \\
            0 & 0 
          \end{matrix} \right )
        = X,
      \end{equation*}
      and such that the restriction of $\psi$ to the diagonal torus of
      $\SL_2$ identifies with the cocharacter $\chi$.
    \item $(\Ad \circ \psi,\glie)$ is a tilting module for $\SL_2$; its
      indecomposable summands are indecomposable tilting modules
      $T(n)$ for $n \le 2p-2$.
    \item $(\Ad \circ \psi,\glie)$ is a semisimple $\SL_2$-module if
      and only if $\Ad \circ \chi$ is a cocharacter of $\GL(\glie)$
      associated with $\ad(X) \in \gl(\glie)$.
    \item If $\ad(X)^{p-1} = 0$, then $\Ad \circ \chi$ is a
      cocharacter of $\GL(\glie)$ associated with $\ad(X) \in
      \gl(\glie)$.
  \end{enumerate}
\end{lem}

\begin{proof}
  The main result of \cite{mcninch-optimal} yields (1).  For (2) see
  \cite{seitz} or \cite{mcninch-optimal}*{Prop. 36}.

  For (3), we first assume $(\Ad \circ \psi,\glie)$ is semisimple.
  Since $T(n)$ is semisimple if and only if $n<p$, $(\Ad \circ
  \psi,\glie)$ is restricted as well.  If we choose a high weight
  vector in each simple summand, it is a consequence of the well-known
  description of restricted semisimple $\SL_2$-modules that this
  collection of vectors is an $\ad(X)$-basis for $\glie$, and that
  $\Ad \circ \chi$ is the cocharacter determined by this
  $\ad(X)$-basis.

  On the other hand, if $(\Ad \circ \psi,\glie)$ is not semisimple,
  then it has an indecomposable summand $T(n)$ for some $p \le n \le
  2p-2$.  Thus the $n$-th weight space of $\chi(\Gm)$ on $T(n)$ is
  non-zero. On the other hand, note that all Jordan blocks of $\ad(X)$
  acting on $\glie$ have size $\le p$. Thus if $\kappa$ is a
  cocharacter of $\GL(\glie)$ associated with $\ad(X)$, then all
  weights $\mu$ of $\kappa(\Gm)$ on $\glie$ satisfy $-p +1 \le \mu \le
  p-1$.  This shows that $\kappa$ and $\Ad \circ \chi$ are not
  conjugate, so that $\Ad \circ \chi$ is not associated to $\ad(X)$.
  This proves (3).

  For (4), note that each non-zero nilpotent element of $\sllie_2$
  acts with partition $(p,p)$ on $T(n)$ for $p \le n \le 2p-2$. Thus
  $\ad(X)^{p-1} = 0$ implies that $(\Ad \circ \psi,\glie)$ is
  semisimple as an $\SL_2$-module so that (4) follows from (3).
\end{proof}

\begin{prop}
  \label{prop:L-prop}
  Assume that the characteristic $p$ of $K$ is 0 or $p > 2h-2$ where
  $h$ is the maximal Coxeter number of a simple component of $G$.
  Then \HypL holds for $G$ using the adjoint representation
  $(\Ad,\glie)$. Moreover, if $p>0$ and if $X \in \glie$ is nilpotent,
  then $Ad(X)^{p-1}=0$.
\end{prop}

\begin{proof}
  Since $p$ is very good for $G$, $\ad:\glie \to \gl(\glie)$ is
  injective.  If $A \in \glie$ is regular nilpotent, and if $\chi$ is
  a cocharacter associated with $A$, then each weight $n$ of
  $\chi(\Gm)$ on $\glie$ satisfies 
  \begin{equation*}
    -2h + 2 \le n \le 2h -2.   
  \end{equation*}
  If $p>0$, our assumption on $p$ means $p-1 \ge 2h-2$; together with
  the condition $A \in \glie(\chi;2)$, it follows that
  $\operatorname{ad}(A)^{p-1}=0$.  Since the regular nilpotent
  elements are dense in the nilpotent variety, one sees that each
  nilpotent element $X \in \glie$ satisfies $\ad(X)^{p-1}=0$.  Part
  (4) of the previous lemma now shows  \HypL to hold for the action
  of $G$ on $(\Ad,\glie)$ as desired.
\end{proof}

\begin{rem}
  In general, the condition in \HypL may fail for the adjoint
  representation.  Indeed, let $X \in \glie$ be regular nilpotent,
  suppose that $X^{[p]} = 0$, and let $\phi:\SL_2 \to G$ be a
  homomorphism determined by $X$ as in (1) of Lemma \ref{lem:optimal}.
  That lemma shows \HypL to fail in case $(\Ad \circ \phi,\glie)$ is
  not semisimple. Semisimplicity fails e.g. in case $G = \SL(n+1)$
  with $p > n > p/2$; indeed, in that case the indecomposable tilting
  $\SL_2$-module $T(2n)$ appears as a summand of $(\Ad \circ
  \phi,\glie)$, and $T(2n)$ is not semisimple since $2n > p$.
\end{rem}

\begin{rem}
  The hypothesis \HypL holds for the symplectic group $\SP(V)$ or the
  special orthogonal group $\SO(V)$ on the natural representation $V$
  provided only that $p=0$ or $p > 2$ (so that $p$ is good for $G$). 
\end{rem}

\begin{rem}
  If $G$ is a group of type $G_2$ and $p=0$ or $p \ge 5$ (so that $p$
  is good for $G$), $\HypL$ holds using the 7 dimensional
  representation $(\rho,V)$ of $G$. In contrast, the condition in
  $\HypL$ holds on the adjoint representation for $G_2$ only when $p >
  2h-2 = 10$.
  
  Note however that if $A \in \glie$ is regular nilpotent, then
  $d\rho(A)$ is regular nilpotent in $\gl(V)$ so that $d\rho(A)^{p-1}$
  only when $p \ge 11$.
\end{rem}

\newcommand\mylabel[1]{#1\hfil}

\end{document}